\documentclass{amsart}
\usepackage{amsmath,amssymb,amsfonts,epsf} 

\newtheorem{Thm}{Theorem} [section]
\newtheorem{prop}[Thm]{Proposition} 
\newtheorem{lem}[Thm]{Lemma} 
\newtheorem{cor}[Thm]{Corollary} 
\newtheorem{defi}[Thm]{Definition} 
\newtheorem{example}[Thm]{Example}
\newtheorem{remark}[Thm]{Remark}

\newcommand{\bthm}{\begin{Thm}}
\newcommand{\ethm}{\end{Thm}}

\newcommand{\bex}{\begin{example}\rm}
\newcommand{\eex}{\end{example}}
\newcommand{\bpr}{\begin{prop}}
\newcommand{\epr}{\end{prop}}
\newcommand{\bco}{\begin{cor}}
\newcommand{\eco}{\end{cor}}
\newcommand{\bde}{\begin{defi}\rm}
\newcommand{\ede}{\end{defi}}
\newcommand{\ble}{\begin{lem}}
\newcommand{\ele}{\end{lem}}
\newcommand{\bre}{\begin{remark}\rm}
\newcommand{\ere}{\end{remark}}

\newcommand{\CC}{\mathcal C}
\newcommand{\M}{\mathcal M}

\newcommand{\N}{\mathbb N}
\newcommand{\R}{\mathbb R}
\newcommand{\Z}{\mathbb Z}

\newcommand{\K}{\mathcal K}

\newcommand{\D}{\mathcal D}
\newcommand{\F}{\mathcal F}

\newcommand{\Om}{\Omega}

\newcommand{\ud}{\underline}

\begin{document}
\author{Paolo Salvatore}
\title{The topological cyclic Deligne conjecture}

\maketitle

\begin{abstract}
Let $O$ be a cyclic topological operad with multiplication. 
In the framework of the cosimplicial machinery by McClure and Smith, 
we prove that the totalization of the 
cosimplicial space associated to $O$
has an action of an operad equivalent to the framed little 2-discs operad.
\end{abstract}

\section{Introduction}

The Deligne conjecture, proved by several authors, states that 
the little 2-discs operad $\D_2$ acts on the chain level on the Hochschild complex of any associative algebra. The question was motivated by a well known action on the homological level  \cite{Ger}.
The natural generalization of the conjecture is about the Hochschild complex of an operad with multiplication, that is an operad with a map from the associative operad. The original statement is recovered taking the endomorphism operad of an associative algebra.
McClure and Smith proved a topological version of the Deligne conjecture \cite{MS99,MS02,MS04}, where the Hochschild complex is replaced by the totalization of a cosimplicial space associated to a topological operad with multiplication. 
A major example, constructed by Sinha \cite{Sin}, occurs when  
the operad with multiplication is a version of the little $n$-discs operad, the Kontsevich operad $\K_n$. 
The (homotopy) totalization of $\K_n$ has as factor the space $Emb_n$ of long knots in $\R^n$ for $n>3$.
We proved in \cite{imrn} that $Emb_n$ itself has an action 
of the little 2-discs operad.
The {\em cyclic} Deligne conjecture states that the operad of framed 
little 2-discs $f\D_2$ acts on the chain level on the Hochschild complex of a Frobenius algebra.
The conjecture was proved in this form by Kaufmann \cite{KaDel}, and Kontsevich-Soibelman \cite{KonSoi}. 
 Since the endomorphism operad of a Frobenius algebra is cyclic, the natural generalization of the conjecture claims an action of the chains of $f\D_2$ on the Hochschild complex of a cyclic operad with multiplication. This statement was proved on the homological level by Menichi \cite{Men}. 
The topological counterpart of the statement says that the framed little 2-discs operad $f\D_2$ acts on the totalization of a topological cyclic operad with multiplication. This is the main result of our paper (theorem \ref{main}). An equivalent statement was conjectured by McClure and Smith \cite{MS04}. Hu \cite{Hu} stated a related result and applied it to the Hochschild complex of commutative Frobenius algebras. 

We will revisit the proof by McClure and Smith of the non-cyclic statement showing that the model of $\D_2$ they used, denoted here $MS$,
is strongly related to Kaufmann's spineless cacti operad. In fact the elements of $MS$ can be thought of as spineless cacti together with a monotone parametrization. As far as we know this has not been observed before. The parametrizations are necessary if we want to include the constant base point at the 0 level of the operad.

We will then construct an operad $fMS$, a model of $f\D_2$, whose elements correspond to Voronov's cacti with a parametrization, and prove the main theorem. 
The proof is conceptually similar to that of the non-cyclic statement, except that the triple defining cyclic operads rather than non-cyclic operads comes into play.
 
The main motivation of this paper is a joint project of 
Ryan Budney and the author to endow the space $fEmb_n$ of framed long knots in $\R^n$ with a meaningful action of the framed little 2-discs operad.  
The action is expected because $fEmb_n$ is the (homotopy) totalization of an operad with multiplication equivalent to the framed little $n$-discs operad  $f\D_n$ \cite{imrn}, and $f\D_n$ is equivalent to a cyclic operad \cite{Bu}.

The paper is relatively self-contained. In sections \ref{two} and \ref{three} we recall the definitions and give several examples of cosimplicial spaces and (non-symmetric) operads. In section \ref{cac} we define geometrically the operad $MS$ and its framed version $fMS$. Then we reprove in section \ref{deligne} the theorem of McClure and\textbf{} Smith (theorem \ref{noncyc}) from this geometric point of view. In section \ref{cyc} we recall the definition of the cyclic category and of a cocyclic space.
In section \ref{ope} we define cyclic operads,
and we prove the main theorem \ref{main}. 
For example the theorem gives an action of $fMS$ on based loop groups.
Another important example we recover is the string topology action of $fMS$
on a free loop space, that makes sense according to Voronov \cite{Vor} in the category of correspondences. 
This paper grew out of my lecture notes at the Lille Summer School on Algebraic Topology in June 2007. I am grateful to the Laboratoire
Painlev\'e and in particular to Sadok Kallel for support and hospitality.

\section{Cosimplicial spaces} \label{two}

We recall the definition of a cosimplicial space.
Let $\Delta$ be the category with the standard ordered sets $[k]=\{0<\dots<k\}$ as objects $(k \in \N)$ 
 and monotone maps as morphisms.
For a fixed $n \in \N$ and $0 \leq i \leq n+1$   let $\delta_i:[n] \to [n+1]$ be the morphism such that 
$\delta_i(j)=j$ if $j<i$ and $\delta_i(j)=j+1$ if $j \geq i$.

For $0 \leq i \leq n$  let $\sigma_i:[n+1] \to [n]$ be the morphism such that
$\sigma_i(j)=j$  if $j\leq i$ and $\sigma_i(j)=j-1$ if 
$j>i$.
It is well known that $\Delta$ is the category generated by these morphisms under the relations
\begin{align*}
\delta_j \delta_i &= \delta_i \delta_{j-1} \quad i<j \\
\sigma_j \sigma_i &= \sigma_i \sigma_{j+1} \quad i\leq j \\
\sigma_j \delta_i &= \delta_i \sigma_{j-1} \quad i<j \\
\sigma_j \delta_i &= id    \quad j \leq i \leq j+1 \\
\sigma_j \delta_i &= \delta_{i-1} \sigma_{j} \quad i>j+1 
\end{align*}
\bde 
A {\em cosimplicial space} is a covariant functor from the category $\Delta$ to the category of topological spaces. The image under such functor of a morphism $\delta_i$ is a {\em coface} denoted $d^i$ and 
the image of a morphism $\sigma_j$ is a {\em codegeneracy} denoted 
$s^j$. 
\ede
The definition is dual to that of a {\em simplicial space}, that 
is a {\em contravariant} functor from $\Delta$ to the category of topological spaces.
\bex \label{delta}
Let us consider the $k$-simplex
$$\Delta^k=\{(x_1,\dots,x_k)\quad|\quad 0\leq x_1\leq \dots \leq x_k \leq 1\}.$$
The collection $\Delta^*$ forms a cosimplicial space with 
the cofaces
\begin{align*}
d^i(x_1,\dots,x_k)&=(x_1,\dots,x_i,x_i,\dots,x_k) \quad 1\leq i \leq k \\
d^0(x_1,\dots,x_k)&=(0,x_1,\dots,x_k) \\
d^{k+1}( x_1,\dots,x_k)&=(x_1,\dots,x_k,1) \\
\end{align*}
and with the codegeneracies
$$s^i(x_1,\dots,x_k)=(x_1,\dots,\widehat{x_{i+1}},\dots,x_k) \quad 0 \leq i \leq k-1\, .$$
\eex

\begin{defi}
The totalization $Tot(S^*)$ of a cosimplicial space $S^*$ is the space of natural transformations of functors
 $\Delta^* \to S^*$. This means that 
 $$Tot(S^*) \subset \prod_{k \in \N} Map(\Delta^k,S^k)$$
 is the space of sequences of maps $f^k:\Delta^k \to S^k$
 commuting with the codegeneracies and the cofaces.
\end{defi}

The definition is dual to that of the {\em realization} of a simplicial space.  
 
We recall two standard examples of cosimplicial spaces.

\bex \label{basedloop}
{\em The cosimplicial model for the based loop space $\Omega X$}

For a given topological space $X$ with base point $* \in X$
let us consider the cosimplicial space 
$\omega X$ such that
$(\omega X)^k = X^k$, with
\begin{align*}
&d^0(x_1,\dots,x_k)=(*,x_1,\dots,x_k) \\
&d^i(x_1,\dots,x_k)=(x_1,\dots,x_i,x_i,\dots x_k) \quad 1 \leq i \leq k \\ 
&d^{k+1}(x_1,\dots,x_k)=(x_1,\dots,x_k,*) \\ 
&s^i(x_1,\dots,x_{k+1})=(x_1,\dots, \widehat{x_{i+1}},\dots,x_{k+1}) \quad 0 \leq i \leq k \,.\\ 
\end{align*}
It is well known that the evaluation maps 
$\Omega X \times \Delta^k \to X^k$  
$$(\gamma, (x_1,\dots,x_k))
\mapsto (\gamma(x_1),\dots,\gamma(x_k))$$
 induce by adjointness  
a homeomorphism $\Omega X \cong Tot(\omega X)$. 
\eex

\bex \label{lx}
{\em The cosimplicial model for the free loop space $LX$}

For a topological space $X$,
let us consider the cosimplicial space $lX$ such that
\begin{align*}
&(lX)^k = X^{k+1}=\{(x_0,\dots,x_k)\,| \,x_i \in X\} \\
&d^i(x_0,\dots,x_k)=(x_0,\dots,x_i,x_i,\dots,x_k)  \quad 0 \leq i \leq k  \\
&d^{k+1}(x_0,\dots,x_k)=(x_0,\dots,x_k,x_0) \\
&s^i(x_0,\dots,x_{k+1})=(x_0,\dots, \widehat{x_{i+1}},\dots,x_{k+1}) 
\quad 0 \leq i \leq k \,.
\end{align*}

The evaluation map $LX \times \Delta^k \to X^{k+1}$  sending
$$(\gamma,(x_1,\dots,x_k)) \mapsto (\gamma(0),\gamma(x_1),\dots,\gamma(x_k))$$
induces by adjointness a homeomorphism $LX \cong Tot(lX)$.

Actually, modulo this homeomorphism, the cosimplicial inclusion $i:\omega X \to lX$ defined by $i^k(x_1,\dots,x_k)=(*,x_1,\dots,x_k)$
corresponds to the inclusion $\Omega X \to LX$.
\eex

\section{Operads} \label{three}

We recall some material on operads.
For more details we refer to \cite{book}.
Roughly speaking an element of an operad is an operation with
many inputs and one output.

\begin{defi} \label{operad}
A non-symmetric topological operad $O$ is a collection of spaces
$O(k), \,k \in \N$, 
 together with
composition maps
$$\_\circ_t\_:O(k)\times O(l) \to O(k+l-1) \quad 1\leq t \leq k$$
such that for
$a \in O(n), \, b \in O(p), \,c \in O(q)\, $ 
\begin{align*}
&(a \circ_i b)\circ_{j+p-1} c = (a \circ_j c)\circ_i b
  & 1 \leq i <j \leq n \\  
&a \circ_i (b \circ_j c) =(a \circ_i b) \circ_{i+j-1} c 
  &1 \leq i \leq n,\, 1 \leq j \leq p \, . 
 \end{align*}
In addition there is a point $\iota \in  O(1)$, the {\em unit}, 
that is a bi-sided unit for the $\circ_t$-operations.
\end{defi} 
 
 The operad is {\em symmetric} if for each natural number $k$ the symmetric group  $\Sigma_k$ acts on $O(k)$ on the right, compatibly with the composition maps  \cite{book}, in the sense that for 
$x \in O(m), y \in O(n), \sigma \in \Sigma_m$ and $\tau \in \Sigma_n$
$$(x \sigma) \circ_i (y \tau) = (x \circ_i y)( \sigma \circ_i \tau)$$
where $\sigma \circ_i \tau \in \Sigma_{m+n-1}$ is the permutation 
exchanging by $\sigma$ $m$ consecutive blocks having all one element except the $i$-th having $n$ elements, and acting by $\tau$ on the $i$-th block.

Here are some examples of topological operads.

\bex \label{endo}
For a topological space $X$, the {\em endomorphism} operad 
$End(X)$ is the symmetric operad defined by $End(X)(n)=Map(X^n,X)$,
such that the composition $\circ_i$ is insertion in the 
$i$-th variable:

for $f \in End(X)(k)$ and $g \in End(X)(l)$
$$f \circ_i g \in End(X)(k+l-1)=Map(X^{k+l-1},X)$$
is
$$(f \circ_i g)(x_1,\dots,x_{k+l-1})=f(x_1,\dots,x_{i-1},
g(x_i,\dots,x_{i+l-1}),x_{i+l},\dots, x_{k+l-1}).$$ 
The identity of $X$ is the unit of $End(X)$.
\eex

We say that a space $A$ is {\em acted on} by a (symmetric) operad $O$, or that it is an
$O$-{\em algebra}, if we are given a morphism of (symmetric) operads
$O \to End(A)$.

The concepts of (symmetric) operads 
and their algebras 
can be defined likewise in any (symmetric) monoidal category.

\bex  \label{coendo}
Dually to example \ref{endo}, for a topological space $X$  
the {\em coendomorphism} operad $Coend(X)$
is the symmetric operad
defined
by 
$$Coend(X)(n)=Map(X,X^n)=Map(X,X)^n \, .$$ 
Given maps  
$f=(f_1,\dots,f_k):X \to X^k$ and $g=(g_1,\dots,g_l):X \to X^l$
the composition is defined by
$$f \circ_i g = (f_1,\dots,f_{i-1}, g_1 \circ f_i ,\dots, g_l \circ f_i,
f_{i+1},\dots,f_k):X \to X^{k+l-1} \, .$$ 
\eex

\bex \label{topmon}
Let $M$ be a topological monoid. Then there is an operad $\underline{M}$ with $\ud{M}(n)=M^n$. The composition is 
$$(x_1,\dots,x_m) \circ_i (y_1,\dots,y_n)=
(x_1,\dots, x_{i-1},x_iy_1, \dots,x_i y_n,x_{i+1},\dots,x_m).$$

Actually the coendomorphism operad of example \ref{coendo}
is the special case with $M=Map(X,X)$ endowed  with the
multiplication $f \cdot g = g \circ f$.

\eex

\bex
The associative operad $Ass$ is the operad such that 
$Ass(n)$ is a point for each $n$. This is a non-symmetric operad
and is the final object in the category of non-symmetric topological operads.
An algebra over $Ass$ is exactly a topological monoid.
\eex

\bex
The {\em little $n$-discs } operad $\D_n$ by Boardman and Vogt is the main 
example that motivated historically the introduction of operads \cite{BV}\cite{May}. 
The space $\D_n(k)$ is the space of $k$-tuples of self-maps of 
the unit $n$-disc $f_i:D_n \to D_n,\, i=1,\dots,k$ obtained by
composing translations and dilatations, i.e.
 $f_i(z)=\lambda_i z + z_i,\, \lambda_i \in \R_+, z_i \in \R^n$, and  
such that $f_i(D_n) \cap f_j(D_n) = \emptyset$ for $i\neq j$.
This is a symmetric suboperad of the coendomorphism operad of the disc $D_n$.
We can associate a picture to an element 
of $\D_n$ drawing the images of the maps.
\eex
$$\epsfbox{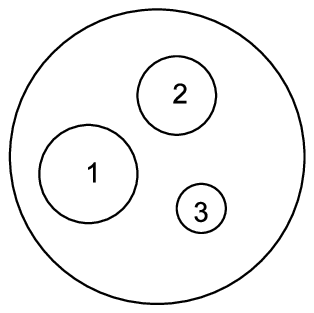}$$ 
The $n$-fold loop space $\Omega^n(X) \cong Map(D_n,S^{n-1};X,x_0)$ 
on a based space $(X,x_0)$ is an algebra over the symmetric operad
$D_n$. The structure map is adjoint to the map
$$\theta_n:D_n(k) \times \Om^n(X)^k \to \Om^n(X)$$
sending $(f_1,\dots,f_k;\omega_1,\dots,\omega_k)$,
for $f_i:D_n \to D_n$ and $\omega_j:D_n \to X$,
to  the $n$-fold loop $$\theta_n(f_1,\dots,f_k;\omega_1,\dots,\omega_k):
D_n \to X$$
coinciding with
 $\omega_i \circ (f_i)^{-1}$ on the image $f_i(D_n)$ for 
 $i=1,\dots,k$  and 
sending any point that is not in the image of any $f_i$
to the base point $x_0 \in X$. 
 
\
A partial converse holds, the 
{\em recognition principle}:
\begin{prop} \cite{BV,May} \label{reco}
Let $A$ be an algebra over $\D_n$, with $n \geq 1$. If $A$ is connected
then $A$ is weakly homotopy equivalent to a $n$-fold loop space.
\end{prop}
More generally this is true if $A$ is {\em group-like}. 
This means that the set of components $\pi_0(A)$,
 that is always a monoid by the action of $\D_n$, is also a group.

\bex
The {\em framed little $n$-discs operad} $f\D_n$ introduced by Getzler 
in \cite{getzler} contains $\D_n$, and is defined similarly as $\D_n$ except that the self-maps in the definition can be composed by oriented rotations.

Precisely $f\D_n(k)$ is the space of $k$-tuples of self-maps
$f_i:D_n \to D_n,\, i=1,\dots,k$  
such that $f_i(z)=\lambda_i A_i z + z_i$ with $\lambda_i \in \R_+$, $A_i \in SO(n)$, 
$z_i \in \R^n$,  and  such that $f_i(D_n) \cap f_j(D_n) = \emptyset$ for $i\neq j$.
Also $f\D_n$ is a suboperad of the coendomorphism operad of the disc $D_n$. There is an obvious inclusion of operads $\D_n \subset f\D_n$. 
\eex

The operad $f\D_n$ acts on the $n$-fold loop space $\Omega^n(X)$ 
of any space $X$ with based $SO(n)$-action.  
\
The recognition principle for algebras over $f\D_n$ is the following.  
\bpr \cite{SW} \label{equireco}
Let $A$ be an algebra over $f\D_n$. If $A$ is group-like then
it is  weakly
 equivalent to $\Omega^n(X)$, where $X$ has a 
based $SO(n)$-action, via a zig-zag of based $SO(n)$-equivariant maps.
\epr

\bde \label{compo}
{\em The language of trees}

Here we follow the terminology of the book by Markl, Shnider and Stasheff \cite{book}, with the important difference that we work
with non-symmetric rather than symmetric operads, and with
planar rather than non-planar trees.
An important point about non-symmetric operads is that it makes sense to compose
elements of the operad along a {\em rooted planar tree}. This is a  planar graph with no cycles, where we allow some open half-edges. 
The {\em root} is one of these half-edges, and the remaining ones are the {\em leaves}. All other edges are called {\em internal} edges.
There is a unique structure of directed graph on a rooted planar tree such that there is a directed path from each point to the root. 
Given a rooted planar tree $T$,
the arity $|v|$ of a vertex $v$ is the number of incoming edges.
The counterclockwise ordering induces an ordering on the incoming edges of a vertex. 
It is customary to draw the tree so that the root is on the bottom and the edges are directed downward.
Here is an example of a tree with 4 vertices of arity
$|x|=0,\; |y|=1,\; |z|=2$ and $|t|=3$.
\epsfxsize=200pt
$$\epsfbox{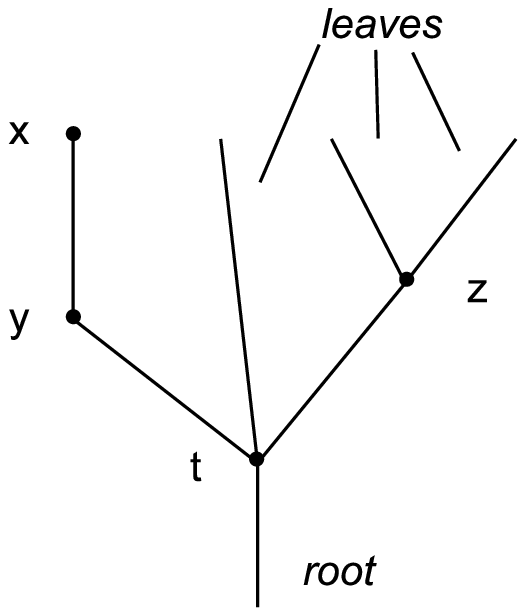}$$
Given a topological operad $O$,
let us consider the product $O(T) = \prod_{v} O(|v|)$ indexed over all vertices $v$ of $T$. This product can be thought of as the space of trees of shape $T$
with vertices labeled by the operad $O$.
If a tree $T'$ is obtained from $T$ by contracting an internal edge
$e:v \to w$, and $e$ is the $i$-th incoming edge of
$w$, then
the composition map $$\circ_i:O(|w|) \times O(|v|) \to O(|w|+|v|-1)$$ induces a map $O(T) \to O(T')$.
The axioms of definition \ref{operad} guarantee that for iterated contractions the result
does not depend on the order. 
In the end there is a well defined composition map
$O(T) \to O(n)$, where $n$ is the number of leaves of the tree $T$.
This can be expressed by the notion of 
a {\em triple}.
We recall that a triple $(F,\mu,u)$ in a category $\mathcal{C}$ is a monoid in the monoidal category of endofunctors of $\mathcal{C}$, i.e. a functor
$F: \mathcal{C} \to \mathcal{C}$ together with natural transformations
$\mu: F \circ F \to F$ and $u: id_{\mathcal{C}} \to F$ satisfying 
associative and unital properties. 
An algebra over the triple is an object $A$ of $\mathcal{C}$ with a morphism $\rho:F(A) \to A$ compatible in an appropriate sense with $\mu$ and $u$.
In our case $\mathcal{C}$ is the category of sequences of topological 
spaces $X=(X(n)) \quad (n \in \N)$ with a distinguished base point $\iota \in X(1)$. The functor $F$ sends a sequence $X$ to 
a sequence $F(X)$,
where $F(X)(n)$ is the space of rooted planar trees on $n$ leaves 
of all possible shapes,
labelled by elements of $X$, modulo the relation deleting arity 1 vertices  labeled by the base point $\iota \in X(1)$. 
The natural transformation $\mu: F \circ F \to F$ sends a tree $T$ with labels that are themselves labeled trees to the labeled tree obtained cutting a neighbourhood of each vertex of $T$ and gluing in its label.  The unit transformation $u:id \to F$ identifies a point to a corolla labeled by the same point, where
a {\em corolla} is a tree with a single vertex.
An algebra over this triple is exactly a non-symmetric topological operad. In other words $F(X)$ is the sequence underlying the free operad generated by $X$. 
\ede

\section{Cacti operads} \label{cac}

We claim that the operads used by various authors dealing with the Deligne conjecture are very similar.

We recall that
two topological operads $A,B$ are said to be
{\em equivalent} if there is a zig-zag of operad
morphisms that are levelwise homotopy equivalences connecting
$A$ and $B$. An $E_n$ operad is an operad equivalent to the little $n$-discs $\D_n$.
Voronov introduced in \cite{Vor} the {\em cacti operad} $Cacti$ that is equivalent to the framed little 2-discs operad $f\D_2$ if we forget constants, i.e. the arity $0$ part. A variation of the cacti operad that is instead equivalent to the little 2-discs $\D_2$, if we forget constants, is the {\em spineless cacti operad}, introduced by Kaufmann
\cite{Ka}.

We will describe in this section an $E_2$ operad constructed by McClure and Smith in \cite{MS02}, that we call $MS$. 
This operad is closely related to 
the 
spineless cacti operad.

The space $MS(n)$ splits as a product of a finite CW-complex
$\F(n)$  and a contractible factor.
We will give a geometric description of $\F(n)$.

\begin{defi} \label{mani}
An element $x \in \F(n)$ is a partition of 
the circle $S^1$ into $n$ closed $1$-manifolds $I_j(x) \;(j=1,\dots,n)$ with equal length (=measure), and pairwise disjoint interiors, such that 
there is no cyclically ordered 4-tuple $(z_1,z_2,z_3,z_4)$ in $S^1$ 
with $z_1,z_3 \in \stackrel{\circ}{I_j}(x),\; z_2,z_4 \in 
\stackrel{\circ}{I_k}(x)$ and $j \neq k$.
\end{defi}

The figure shows an element of $\F(4)$.
$$\epsfbox{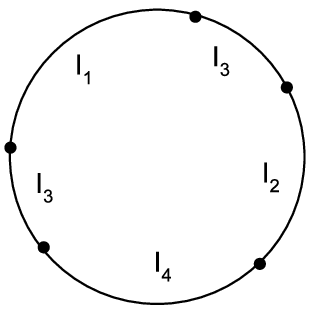}$$
We can define on $\F(n)$ a topology induced by a metric such that two partitions $x,y \in \F(n)$ are close if the manifolds $I_j(x)$ and $I_j(y)$ are close for 
each $j$. For example the metric
$d(x,y)=\sum_{j=1}^n l(I_j(x)\setminus \stackrel{\circ}{I_j}(y) )$ 
does the job, where $l$ denotes the measure.
The symmetric group $\Sigma_n$ acts on $\F(n)$ by reindexing the 
labels.

Clearly $\F(1)$ is a point.
Let $\pi:I \to I/\partial I$ be the projection.
The space $\F(2)$ is $\Sigma_2$-equivariantly homeomorphic to a circle 
with antipodal $\Sigma_2$-action. 
The homeomorphism $I/\partial  I \cong \F(2)$ sends
$\pi(x)$ into the pair of intervals $\{\pi[x,x+\frac{1}{2}],\pi[x+\frac{1}{2},x+1]\}$.

\medskip

{\em The cell structure}

The pullback of a partition $x \in \F(n)$  
to $I$ via $\pi:I \to I/\partial I \cong S^1$ determines
a partition of $I$ into say $k$ intervals labelled by natural numbers between $1$ and $n$. Any component of $\pi^{-1}(I_j(x))$ is labelled
by $j$. 
This gives an ordered sequence of boundary points 
$0=x_0<x_1<\dots<x_{k-1}<x_k=1$, together
with the sequence of labels $X=(X_1, \dots ,X_k)$, such that
$X_i$ is a label of the interval $[x_{i-1},x_i]$.

It follows that $\F(n)$ is a regular CW-complex with one cell
for each fixed sequence of labels $X$. The closed cell indexed by 
a sequence $X$ of length $k$ corresponds to
the subspace of the $(k-1)$-simplex containing the elements    
$$0=x_0 \leq x_1 \leq \dots \leq x_{k-1} \leq x_k =1$$
such that
$$\sum_{i| X_i=j}(x_i-x_{i-1})=1/n$$ for each $j=1,\dots,n$. This subspace is 
homeomorphic to the product of simplices $\prod_{j=1}^n \Delta^{d(j)}$,
where $d(j)$ is one less than the number of components (intervals) of $\pi^{-1}(I_j(x))$.
\
The upshot is that cells of $\F(n)$ correspond to 
those finite sequences $X=(X_i)$ with values in $\{1,\dots,n\}$ such that

\begin{enumerate}

\item All values between 1 and $n$ appear in $X$ 
\item Two adjacent values in $X$ are distinct
\item The sequence $X$ does not have a subsequence of the form $\{i,j,i,j\}$ with $i \neq j$.

\end{enumerate}

Given a partition $x \in \F(n)$ and an index $j$, let us consider
the quotient space of $S^1$ under the relation that identifies points in the same component (closed interval) of $\overline{S^1 \setminus  I_j(x)}$. The result is a circle with a base point, the image of the base point of $S^1$. If we rescale this circle by the factor $n$ and identify it to $S^1$, then we see the quotient map as a based self-map $\pi_j:S^1 \to S^1$.
 
We say that  $c(x):=(\pi_1,\dots,\pi_n):S^1 \to (S^1)^n$ is
the {\em cactus map} and its image
is the {\em cactus} $C_x$ associated to $x$. 

The map $c: \F(n) \to Map(S^1,(S^1)^n)$ is an embedding, but
the image does not form a suboperad of the coendomorphism operad $Coend(S^1)$ of example \ref{coendo}.

The cactus $C_x$ is the union of $n$ circles, the {\em lobes} $S_j(x) \, (j=1,\dots,n)$, that are the images under $c(x)$ of the submanifolds
$I_j(x)$.
The base point $P \in C_x$ is the image of the base point $1 \in S^1$ under the cactus map $c(x):S^1 \to C_x$.

Here is an example of a cactus with the indexes of the lobes.
\epsfysize=100pt
$$\epsfbox{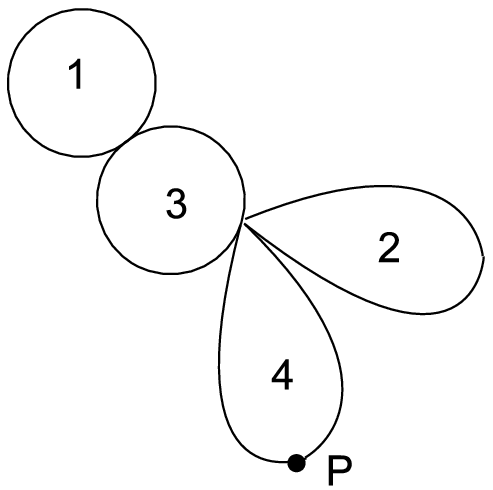}$$
\bde  \label{tree}
{\em The associated tree}

One can associate to each partition $x \in \F(n)$ a rooted planar tree (constant on open cells) constructed as follows. The tree is the graph that has 
a vertex for each lobe (i.e. for each number between 1 and $n$), a vertex for each
intersection point of distinct lobes, and an edge between the vertex of an intersection point and the vertex of each lobe it belongs to.
Let $P \in C_x$ be the base point.
If $P$ is an intersection point then the {\em root vertex} is its vertex. Otherwise $P$ belongs to a unique lobe and 
the root vertex is the vertex of that lobe.
The {\em root} is an additional half-open edge incident to 
the root vertex.
The tree is a ribbon graph with a cyclic ordering of the edges 
incident to a given vertex. The ordering is induced by the standard 
orientation of the circle via the map $c(x):S^1 \to C_x$.
The tree admits a planar embedding such that the cyclic ordering is compatible with the standard orientation of the plane\cite{Vor}.
It is possible to embed the cactus $C_x$ into the plane, in a way that this is compatible with the embedding of the tree on the intersection points, as shown in the next figure. 
The tree is the union of corollas
inscribed in the lobes.
\ede
 
\epsfysize=100pt
$$\epsfbox{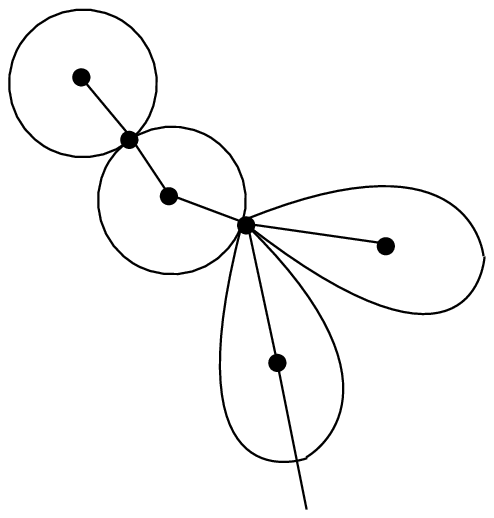}$$

\medskip

{\em More operads}
 
McClure and Smith constructed in \cite{MS99} first an operad
$\CC'$ that is equivalent to little 2-discs $\D_2$ away from  
the arity $0$ part, because $\CC'(0)$ is empty but $\D_2(0)$ is a point.
As a $\Sigma_n$-space $\CC'(n)=\F(n) \times \stackrel{\circ}{\Delta}_{n-1}$.
 
The space $\CC'(n)$ is naturally homeomorphic to a space of partitions
defined similarly as in definition \ref{mani}, except that one drops the condition that all manifold must have the same length.
In this setting there is also a cactus map $\CC' \to Coend(S^1)$ 
that is an embedding onto a symmetric suboperad of $Coend(S^1)$ and defines 
an operad structure on $\CC'$.

The spineless cacti operad $Cact$ by Kaufmann \cite{Ka}
is actually a pseudo-operad, in the sense that is associative,
non-unital,
and is isomorphic to $\CC' \times R$, where $R$ is the 
constant pseudo-operad $R(n)=\R_+$ with composition $\R_+ \times
\R_+ \to \R_+$ equal to the projection onto the first factor. 
The extra parameter corresponds to the length of the parametrizing
circle in \cite{Ka}.

\medskip

McClure and Smith extended in \cite{MS99} 	the operad $\CC'$, that has 
$\CC'(0)$ empty, to an operad
$\CC$ equivalent to the little $2$-discs operad such that 
$\CC(0)$ is a base point. The operad $\CC$ is the smallest suboperad
of $Coend(S^1)$ containing $\CC'$ and the base point in arity 0.

Later on McClure and Smith defined a larger operad that is
denoted $\D_2$ in \cite{MS02}. We present a geometric construction 
of such operad. 

\begin{defi}
Let $Mon(I,\partial I)$ be the space of self-maps of the circle
$S^1 \cong I/\partial I$ induced by weakly monotone self-maps of the interval $I$ restricting to the identity on the boundary $\partial I$. 
We define $MS(n)=\F(n) \times Mon(I,\partial I)$ for $n>0$.  
We also define $MS(0)$ to be a point.
\end{defi}

\bre
The space $Mon(I,\partial I)$ is homeomorphic to the contractible space $Tot(\Delta^*)$, i.e. the monoid of 
cosimplicial endomorphisms of the cosimplicial space $\Delta^*$.
The homeomorphism $T:Mon(I,\partial I) \to Tot(\Delta^*)$ 
sends $f:I \to I$ to the sequence of maps $T(f)^k:\Delta^k \to
\Delta^k$ with 
$$T(f)^k: (0 \leq x_1 \leq \dots \leq x_k \leq 1) \mapsto (0 \leq f(x_1) \leq \dots \leq f(x_k) \leq 1).$$
\ere 

There is an embedding $\zeta_n:MS(n) \to Coend(S^1)(n)$  sending, for $n>0$, 
$x \in \F(n)$ and $f \in Mon(I,\partial I)$ to the composition
$c(x) \circ f:S^1 \to (S^1)^n$. 
The following characterization is immediate.
\begin{prop} \label{conditions}
The image $\zeta_n(MS(n))$ is the space 
of {\em based} maps $g:S^1 \to (S^1)^n$ such that
\begin{enumerate}  
\item $\pi_i \circ g$
is a weakly monotone degree 1 map for each $i=1,\dots,n$
\item there is a partition of $S^1$ into closed intervals intersecting at the boundary
such that on each interval $\pi_i \circ g$ is constant  
for all indexes $i$ except one special index;
\item The counterclockwise sequence of special indices does not contain
a subsequence of the form $\{i,j,i,j\}$ with $i \neq j$;  
\end{enumerate}
\end{prop}
It follows that the image of the embedding $\zeta$ is a symmetric suboperad of $Coend(S^1)$. This gives 
an operad structure on $MS$.
\

\begin{prop} \label{msd}
The operad $MS$ is isomorphic to the $E_2$ operad denoted 
$\D_2$ in \cite{MS02}.
\end{prop}

A proof will be presented in the appendix.

\

We build next an operad $fMS$ that is equivalent to the framed little
2-discs. It is related to $MS$ in the same way as Voronov's cacti operad is related to Kaufmann's spineless cacti operad.
 
\bde \label{fms}
Let us consider the $\Sigma_n$-space $fMS(n) = MS(n) \times (S^1)^n$ for any natural number $n$.
There is an embedding $c':fMS(n) \to Coend(S^1)(n)$ 
sending 
$$(x,f,(z_1,\dots,z_n)) \mapsto (L_{z_1} \times \dots \times 
L_{z_n}) \circ c(x) \circ f$$
where $L_z:S^1 \to S^1$ is multiplication by $z \in S^1$.
A map $g:S^1 \to (S^1)^n$ belongs 
to the image $c'(fMS(n))$ if and only if it is an {\em unbased} map
satisfying the conditions 1,2,3 of proposition \ref{conditions}.
It follows that the image of $c'$ is a symmetric suboperad of $Coend(S^1)$, 
and this defines an operad structure on 
$fMS$. 
\ede

In order to prove that $fMS$ is equivalent to $f\D_2$ we need 
a lemma. The lemma is exposed as Theorem 2.2.2 in \cite{CohVor}
and follows from Theorem 7.3 in \cite{SW}, that is a generalization
of a similar result by Fiedorowicz for the little 2-discs operad.

\ble \label{ribbon}
Let $O$ be a symmetric operad together with
a $\Sigma_k$-equivariant  
homotopy equivalence $O(k) \simeq f\D_2(k)$ for each $k$,  and
together with an operad map $E \to O(k)$ from an $E_1$-operad $E$.
Then $O$ is equivalent to $f\D_2$.
\ele

\bpr
The operad $fMS$ is equivalent to the framed little 2-discs operad
$f\D_2$.
\epr

\begin{proof}
We know from Proposition \ref{msd} that $MS$ is an $E_2$-operad, and this gives for each $k$
a $\Sigma_k$-equivariant homotopy equivalence
$$fMS(k)=MS(k) \times (S^1)^k \simeq \D_2(k) \times (S^1)^k 
= f\D_2(k).$$
In addition $MS \subset fMS$ contains an $E_1$-operad $MS_1$ 
defined in the appendix.
By lemma \ref{ribbon} this concludes the proof.
\end{proof}

For $y=(x,f,(z_1,\dots,z_n)) \in fMS(n)$ 
we denote by $C_y$ and call {\em cactus} the image of 
$(L_{z_1} \times \dots \times L_{z_n})\circ c(x):S^1 \to (S^1)^n$.
We also call {\em lobe} and denote $S_j(y)$ the image of $I_{j}(x)$ via the same map. The lobe has a standard identification to $S^1$ via the
projection $(S^1)^n \to S^1$ onto the $j$-th coordinate. The base point of the cactus is the image of the base point of $S^1$.
In general the base points of the lobes are not intersection points or the base point of the cactus as in the spineless case.

The connection to Voronov's cacti operad $Cacti$, that actually is a
pseudo-operad, is as follows.
The operad $fMS$ contains an operad $f\CC'$ with $f\CC'(n)=\CC'(n) \times (S^1)^n$. In particular $f\CC'(0)$ is empty.
 The inclusion $f\CC'(n) \to fMS(n)$ is a 
 $\Sigma_n$-equivariant homotopy equivalence
for all $n>0$, and there is an isomorphism of pseudo-operads
$Cacti \cong f\CC' \times R$.

\medskip

\section{The non-cyclic Deligne conjecture} \label{deligne}

We shall reprove in this section the result by McClure and Smith.
We need the following definition.
\begin{defi}
A topological operad with multiplication $O$ is
a non-symmetric operad $O$ together with a morphism 
$Ass \to O$ from the associative operad. 
\end{defi}

Equivalently $O$ is a non-symmetric operad in 
the monoidal category of based spaces with the cartesian product as tensor product.
Also equivalently $O$ is a topological operad
with points $m \in O(2)$ (the multiplicaton)
and $u \in O(0)$ (the unit of the multiplication) such that
\begin{align*}
m \circ_1 m &= m \circ_2 m \in O(3) \\
m \circ_1 u &= m \circ_2 u = \iota \in O(1) 
\end{align*}

\bre \label{pointedtriple}
There is a triple $F_*$ in the category of sequences of based spaces
such that an operad with multiplication is exactly an algebra over 
$F_* \,$. The triple $F_*$ is closely related to the triple $F$ of definition 
\ref{compo}. Namely, for a sequence $X$, $F_*(X)$ is
the quotient of $F(X)$ under the identification relation induced by collapsing internal edges with both vertices labeled by base points.
\ere

We can relate now cosimplicial spaces and operads.

\begin{defi}
An operad $O$ with multiplication defines a cosimplicial space $O^*$ sending $[k]$ to $O(k)$.
The coface operator $d^i:O(k) \to O(k+1)$ is defined by  
\begin{align*}
&d^i(x)=x \circ_i m\quad \text{for $1\leq i \leq k$}  \\
&d^0(x)=m \circ_2 x \\
&d^{k+1}(x)=m \circ_1 x
\end{align*} 
The codegeneracies $s^i:O(k)\to O(k-1)$  are defined by
$$s^i(x)=x \circ_{i+1} u $$ 
\end{defi}

\begin{Thm} {\rm (McClure-Smith)} \label{noncyc} \cite{MS02}
Let $O$ be an operad with multiplication.
Then the totalization $Tot(O^*)$ 
 admits the
 action of the $E_2$ operad 
   $MS$. 
   \end{Thm} 
  
\begin{proof}
   
We explain geometrically the action of $MS$ on $Tot(O^*)$.

We need to describe
the action of $MS$ on $Tot(O^*)$ via maps
$$\theta_m:MS(m) \times Tot(O^*)^m \to Tot(O^*).$$

Given $(x,f) \in MS(m)=\F(m) \times Tot(\Delta^*)$ and $m$ cosimplicial maps 
$$a_1^*,\dots,a_m^* \in Tot(O^*)$$ we want to define 
a cosimplicial map 
$$a^* =\theta_m(x,f;a_1^*,\dots,a_m^*) \in Tot(O^*),$$
 so we need compatible maps $a^k:\Delta^k \to O(k)$ for every non-negative integer $k$. A point of the simplex $\Delta^k$ 
is a sequence of points of the unit interval 
$$0=x_0 \leq x_1\leq \dots \leq x_k \leq 1 \,.$$

The cactus map induces a curve $$\gamma = c(x) \circ f \circ \pi: I \to C_x \subset (S^1)^m \, .$$

We construct a rooted planar tree $T$ associated to these data.
We start with the tree associated to $x \in \F(m)$
in definition \ref{tree} with its root removed. This tree has 
\begin{itemize}
\item a vertex for each lobe
\item a vertex for each
intersection point
\item edges between vertices of intersection points
and vertices of their lobes 
\end{itemize}

We call the points of the form
$\gamma(x_j)\quad (0\leq j \leq k)$ {\em special points}.
For each special point that is not an intersection point we add a new vertex, and an edge connecting it to the vertex of the lobe containing $\gamma(x_j)$.
 
Finally for each point of the form $x_j \quad (0\leq j \leq k)$ we add a leaf (half-open edge) incident to the vertex of $\gamma(x_j)$.
The edge incident to the vertex of the base point $\gamma(x_0)$ is
the root. The difference with the tree in definition \ref{tree} is that we
subdivide the former root inserting a vertex. 

The tree $T$ is a planar rooted tree: the ordering of the edges incident to a vertex is induced by the standard orientation of the interval via $\gamma$ and by the ordering of the indices. Both the tree $T$ and the cactus $C_x$ can be 
embedded together consistently into the plane, so that the embeddings are compatible on points of the form $\gamma(x_j)$ and on intersection points.
An example with $k=6$ and $m=4$ is given in the next figure.
$$\epsfbox{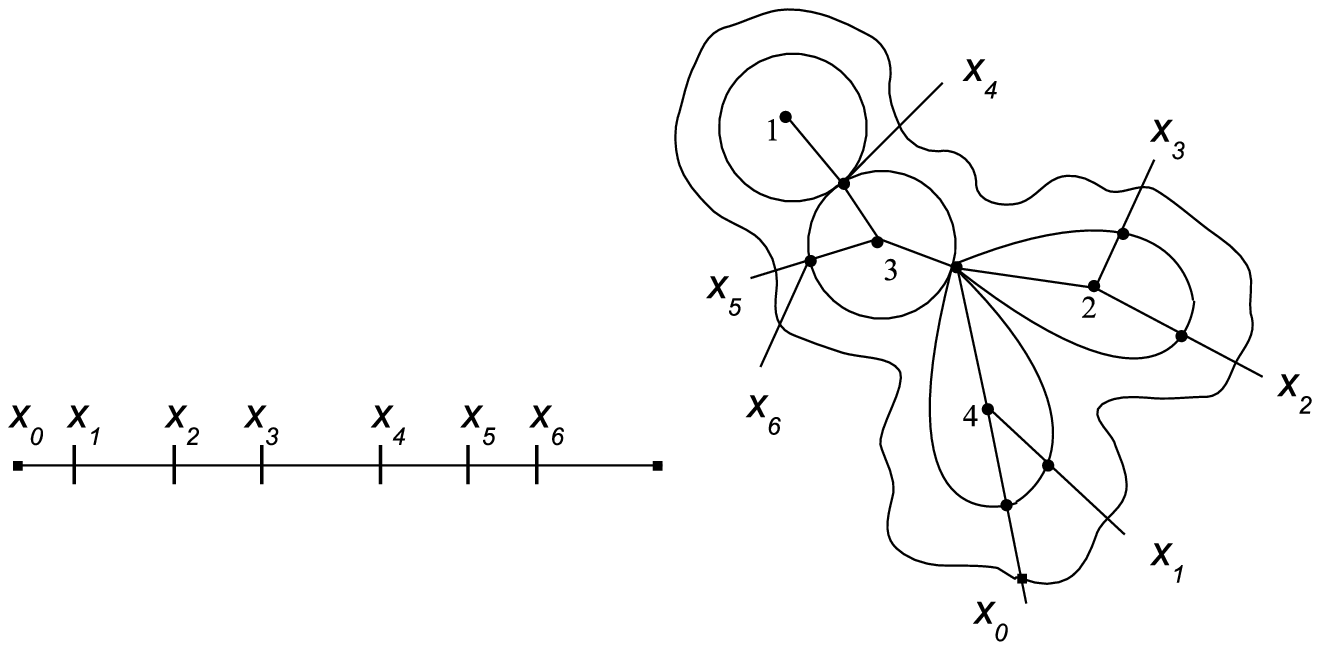}$$
Next we label each vertex $v$ of $T$ by an element of the operad
of the same arity.

If $v$ corresponds to an intersection point of lobes, or to a point 
of the form $\gamma(x_j)$, then we label it by the base point of $O(|v|)$.
Let $Int \subset C_x$ be the set of all intersection points.

If $v$ corresponds to a lobe $S_j$, then
 let us take the set 
 $$H_j = (Int \cup \gamma (\{x_1,\dots, x_k \})) \cap (S_j-*)$$
 of special points and intersection points in the lobe that are not its base point. 

Suppose that $H_j$ has $l_j$ points.
The arity of $v$ is exactly $|v|=l_j$.
There is a standard identification of the lobe $S_j$ with $S^1$
via inclusion $S_j(x) \subset C_x \subset  (S^1)^m$ followed by projection onto the $j$-th factor.
Pulling back $H_j$ via $I \to I/\partial I \cong S^1 \cong S_j$ gives 
$l_j$ points in the unit interval $I$, i.e.
a point $y_j$  of the simplex 
$\Delta^{l_j}$. Actually $y_i$ is in the interior of the simplex. The cosimplicial map $a_j \in Tot(O^*)$ in degree $l_j\;$
$a_j^{l_j}:\Delta^{l_j} \to O(l_j)$ evaluated on $y_j$
gives an element of the operad $a_j^{l_j}(y_j) \in O(l_j)$ that we assign as label of $v$.
In the terminology of definition \ref{compo}
our labelled tree is some $\Gamma \in O(T)$.
The generalised composition $O(T) \to O(k)$ applied to
the labelled tree $\Gamma$ defines the desired element
$a^k(x_1,\dots,x_k) \in O(k)$. 

\medskip

The resulting map $a^*:\Delta^* \to O^*$
is a cosimplicial map essentially
because $a_1^*,..,a_m^*$ are all cosimplicial maps. 
The maps $\theta_m$ define a $MS$-algebra structure on $Tot(O^*)$. 
To see this we need to look at the labelled trees
whose composition defines the structure maps, and use the fact that the ordering of composition by contracting edges is not relevant.
More precisely we are using the fact that $O$ is an algebra over the 
triple $F_*$ of remark \ref{pointedtriple}.
\end{proof}

\bex
Suppose that $x \in MS(2)$ , so that $C_x$ has two lobes,
that the base point is in $S_1$, and that $k=1$.
We wish to compute the value of the map $a^1:\Delta^1 \to O(1)$ on $x_1 \in I$.  
We suppose that the intersection point $Q =S_1 \cap S_2$ coincides 
with the special point $Q=\gamma(x_1)$. 
We also assume that the curve $\gamma$ at $x_1$ goes from the lobe $S_1$ into $S_2$, i.e. $\gamma([0,x_1]) \subset S_1$.   
Then $l_1=1, \, l_2=0$ and $y_1$ is the point mapping to $Q$ via
$I \to S^1 \cong S_1$.
It follows that $$a^1(x_1)=(a_1^1(y_1) \circ_1 m_2) \circ_2 a_2^0(*).$$

\epsfxsize=150pt
$$\epsfbox{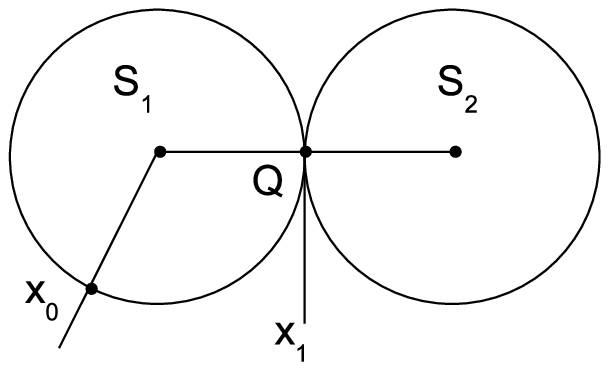}$$

\eex

Theorem \ref{noncyc} and the recognition principle (proposition \ref{reco})
immediately imply:

\bco \label{double}
If $O$ is an operad with multiplication and $Tot(O^*)$ is group-like then $Tot(O^*)$ has the weak homotopy type of a double loop space.
\eco

\bex
Let $M$ be a topological monoid. The operad $\underline{M}$ 
of example \ref{topmon}
has a multiplication by choosing $m_2=(1,1) \in M^2$, and
$m_0=* \in M^0$.
Theorem \ref{noncyc} tells us that the operad $MS$ acts
on $Tot(\omega M) \cong \Omega M \simeq \Omega^2(BM)$.
This is consistent with corollary \ref{double}. 
\eex

\bex
{\em Knots from operads}

We recall that a {\em long knot} is an embedding 
$\R \to \R^n$ that is fixed near infinity, 
and a {\em framed long knot} is a long knot with 
a trivialization of its tubular neighbourhood that is fixed near infinity.
The Kontsevich operad $\K_n$ is an operad
with multiplication that is equivalent to the little discs
$\D_n$, which instead does not admit a multiplication.
Sinha proved in \cite{Sin}, using Goodwillie calculus, that the {\em homotopy} totalization $\widetilde{Tot}(\K_n^*)$ is
weakly equivalent to $Emb_n \times \Omega^2 S^{n-1}$ for $n>3$, where $Emb_n$ is the 
space of long knots in $\R^n$.
The homotopy totalization of a cosimplicial space is weakly equivalent to the homotopy limit of the cosimplicial space as a diagram.
One cannot apply directly theorem \ref{noncyc} to
the homotopy totalization, but
there is a derived version of theorem \ref{noncyc} by McClure and Smith (theorem 15.3 in \cite{MS04}) constructing an $E_2$ operad that acts on the homotopy totalization of any operad with multiplication.
By corollary \ref{double} the space $Emb_n \times \Omega^2 S^{n-1}$,
which is connected by Whitney's theorem,
is a double loop space for $n>3$.
We proved in \cite{imrn} that $Emb_n$ itself is a double loop space for $n>3$. We also proved that the space of framed long knots $fEmb_n$ 
is a double loop space,
because it is weakly equivalent to $ \widetilde{Tot}(f\K_n^*)$, where $f\K_n$ is the {\em framed} Kontsevich operad, that is equivalent to the framed little $n$-discs $f\D_n$.

\eex

\section{The cyclic category} \label{cyc}

The cyclic category $\Lambda$ was introduced by Connes, and is an 
extension of the category $\Delta$ by adding for each $n$ an automorphism $\tau_n:[n] \to [n]$ of order $n+1$, such that
\begin{align*}
&\tau_n \delta_i = \delta_{i-1} \tau_{n-1} \quad     1 \leq i \leq n \\
&\tau_n \sigma_i = \sigma_{i-1} \tau_{n+1} \quad     1 \leq i \leq n 
\end{align*}
As a consequence of the axioms $\tau_n \delta_0 = \delta_n$
and $\tau_n \sigma_0 = \sigma_n (\tau_{n+1})^2 $.  
We refer to \cite{Drin} for a discussion of the cyclic category.

A {\em cocyclic} space is a covariant functor from $\Lambda$ to the category of topological spaces. In particular by restriction to $\Delta \subset \Lambda$ a cocyclic space has an underlying 
cosimplicial space, that we denote with the same name.
We denote by $t_n$ the image of $\tau_n$ under the functor
defining a cocyclic space.

John Jones proved the following result.

\begin{prop} \cite{JJ}\label{jones}
Let $C$ be a cocyclic space.
Then the space $Tot(C)$ has a natural circle action.
\end{prop}

\bex \rm \cite{JJ} \label{freeloop}
The cosimplicial model of the free loop space $lX$ 
of example \ref{lx} has
a compatible cocyclic space structure 
defined by maps $t_n:X^{n+1} \to X^{n+1}$ with $$t_n(x_0,\dots,x_n)=(x_1,\dots,x_n,x_0).$$
The resulting circle action on $Tot(lX)$,
modulo the identification $Tot(lX) \cong LX$, 
coincides 
with the reparameterization action on the free loop space  
$S^1 \times LX \to LX$ defined by $(z , \gamma) \mapsto 
\gamma \circ L_z$, where $L_z:S^1 \to S^1$ is multiplication by $z$. 
\eex

\bex
Let $G$ be a topological group.
The cosimplicial space $\omega G$ of example \ref{basedloop} has a compatible cocyclic space structure defined by maps $t_n:G^n \to G^n$ with
$$t_n(g_1,\dots,g_n)= (g_1^{-1}g_2,\dots,g_1^{-1}g_n,g_1^{-1}).$$
The homeomorphism $\Omega G \cong Tot(\omega G)$ is $S^1$-
equivariant, where the circle action
$\theta:S^1 \times \Omega G \to \Omega G$
 on $\Omega G=Map_*(S^1,G)$  
$$\theta(z, \gamma)(w) = \gamma(z)^{-1} \gamma(zw)$$
is given  by normalized reparameterization.
\eex

\bex
The cocyclic space $\Lambda^*$ is defined so that
$\Lambda^n \subset(S^1)^n$ is the subspace 
of those $(n+1)$-tuples $(x_0,\dots,x_n)$ that are cyclically ordered with respect to the counterclockwise cyclic ordering. 

Similarly as for $\Delta^*$ (example \ref{delta} ) 
the $i$-th coface doubles the coordinate $x_i$, 
the $i$-th codegeneracy removes the coordinate $x_{i+1}$
and the automorphism $t_n$ acts cyclically on the coordinates by
$$t_n(x_0,\dots,x_n)=(x_1,\dots,x_n,x_0)$$. 
There is a homeomorphism $\Delta^n \times S^1 \cong \Lambda^n$
sending 
$$((x_1,\dots,x_n) , z) \mapsto (z,z+x_1,\dots,z+x_n)$$
for $z \in \R/\Z \cong [0,1]/\{0,1\} \cong S^1$.
The inclusion $\Delta^n \times \{1\} \to \Lambda^n$ makes 
$\Delta^*$ into a cosimplicial subspace of $\Lambda^*$. 

\eex

The totalization $Tot(X^*)$ of a cocyclic space
$X^*$ is naturally isomorphic to the space $Mor(\Lambda^*,X^*)$
of natural transformations
$\Lambda^* \to X^*$  of cocyclic spaces. 
In particular for any space $Y$ there is a homeomorphism $LY \cong  Mor(\Lambda^* \to lY)$
sending a loop $\gamma$ to the sequence of maps $(x_0,\dots,x_k) \mapsto (\gamma(x_0), \dots,\gamma(x_k))$.
Jones' circle action on $Tot(X^*)$ of proposition \ref{jones}
is induced by the circle action
on $\Lambda^* \cong S^1 \times \Delta^*$.

\section{Cyclic operads and the Deligne conjecture} \label{ope}

Roughly speaking a cyclic operad exchanges the roles of inputs and outputs. 
We give the definition by default in the non-symmetric sense.

\begin{defi} \label{ciclico}
A cyclic operad $O$ is an
operad such that $O(n)$ has an action of the cyclic group
$\Z_{n+1}$ for each $n$, with the generator acting by
$t_n:O(n) \to O(n)$, 
such that
for $f \in O(m), \, g \in O(n)$
\begin{align*}
t_{m+n-1}(f \circ_1 g)&=t_n(g) \circ_n t_m(f) \\
t_{m+n-1}(f \circ_i g)&= t_m(f) \circ_{i-1} g   \quad\quad 1<1\leq m \\  
t_1(\iota)&=\iota
\end{align*}
\end{defi}

A cyclic {\em symmetric} operad is a cyclic operad such that $O(n)$ admits the action of the symmetric group $\Sigma_{n+1}$, extending the action of
$\Z_{n+1}$ and of $\Sigma_n$.

The simplest example of a cyclic operad is the associative
operad $Ass$ with trivial action of the cyclic groups.

\bex
A non trivial example of cyclic non-symmetric operad is the $A_\infty$  operad by Stasheff
\cite{book}. In particular $\Z_5$ acts on the pentagon $A_\infty(4)$
by rotations. 
\eex
\bex
An example of symmetric cyclic operad is the operad $\M$ of moduli space of Riemann surfaces with boundary.
Here $\M(n)$ is the moduli space of the surfaces with $n+1$ boundary
components. The operadic
composition goes by gluing surfaces at the boundary components, and  $\Sigma_{n+1}$ acts on the 
labels of the components.
The genus zero part $\M_0$ is a cyclic operad equivalent to the 
framed little discs operad $f\D_2$.
\eex

\bex
Budney defines in \cite{Bu} the operad of {\em conformal balls} $\mathcal{CB}_n$
that is symmetric cyclic, and equivalent to the framed little $n$-discs
$f\D_n \,$.

\eex

\bde \label{ciclo}
{\em The triple for cyclic operads} 

Getzler and Kapranov described in \cite{gk} the triple for symmetric cyclic operads. We describe the analog triple in the non-symmetric case.
A {\em cyclic sequence} is a sequence of spaces $X(n),\, n \in \N$
together with an action of $\Z_{n+1}$ on $X(n)$ for each $n$. We also assume that $X(1)$ has a base point in $X(1)$ fixed by the $\Z_2$-action.
There is a triple $F^+$ on the category of cyclic sequences. Given a cyclic sequence $X$, if we forget the actions of the cyclic groups then $F^+(X)$ agrees with $F(X)$. By definition $F(X)(n)$ is the space of rooted labelled trees with $n$ leaves, or equivalently with $n+1$ half edges.
The group action of $\Z_{n+1}$ on $F(X)(n)$ 
switches the roots. The action of $[m] \in \Z_{n+1}$ on some
$T \in F(X)(n)$
sets the new root as the $m$-th half-edge in the counterclockwise cyclic 
ordering after the old root. The new root defines a new structure of directed graph on the labelled tree $T$. For each vertex $v$ of arity $p$, or equivalently of valence $p+1$, there is a unique outgoing edge $v_o$ (respectively $v_n$) in the old (respectively new) directed graph structure. Suppose that $v_n$ 
is the $q$-th edge after $v_o$ in the counterclockwise ordering of the 
edges incident to $v$. 
The new label of $v$ is the result of the action of $[q] \in \Z_{p+1}$ 
on the old label in $X(p)$. 
Given a cyclic operad $O'$, by virtue of definition \ref{ciclico}, 
the action of the triple $F(O') \to O'$ respects the action of the cyclic groups. This shows that $O'$ is an algebra over the triple $F^+$.
Viceversa an algebra $O$ over the triple $F^+$ is a cyclic operad.
The definition \ref{ciclico} is satisfied by restricting 
the action $F(O) \to O$ to labelled trees with one internal edge. 
\ede

\bde
A {\em cyclic operad with multiplication} is an
operad $O$ together with a morphism of cyclic operads 
$Ass \to O$.
\ede
Equivalently $O$ is an operad with multiplication $m_2 \in O(2)$, and
a cyclic operad, such that $t_2 (m_2) = m_2$.

\medskip

There is a triple $F^+_*$ describing cyclic operads with multiplication.
The triple $F^+_*$ is defined 
on the category of cyclic sequences of based spaces, with the base points fixed by the cyclic group action. 
If we ignore the cyclic group action, then $F^+_*$ agrees with the triple $F_*$ of remark \ref{pointedtriple}.
The action of the cyclic groups is defined similarly as for $F^+$ in definition \ref{ciclo}.

We are now ready to state the main theorem.

\bthm \label{main}
Let $O$ be a cyclic operad with multiplication.
Then the operad $fMS$ acts on $Tot(O^*)$.
\ethm

\begin{proof}

We need to describe
the action of $fMS$ on $Tot(O^*)$ via maps
$$\theta_m:fMS(m) \times Tot(O^*)^m \to Tot(O^*).$$
Given $x \in fMS(m)$, and cocyclic maps 
$$a_1^*,\dots,a_m^* \in Tot(O^*)=Mor(\Lambda^*,O^*),$$ we want to define a cocyclic map
$$a^*:=\theta_m(x;a_1^*,\dots,a_m^*) \in Tot(O^*)=Mor(\Lambda^*,O^*)$$
 so we need compatible maps $a^k:\Lambda^k \to O(k)$ for $k \in \N$. 

The construction is similar to the one for non-cyclic operads
in the proof of theorem \ref{noncyc}.

Let $\gamma=c'(x):S^1 \to (S^1)^m$ be the parametrization of the cactus 
$C_x$ associated to $x$ in definition \ref{fms}. 
Given a system of cyclically ordered points $(x_0,\dots,x_k)\in \Lambda^k$
on the circle $S^1$ we obtain via $\gamma:S^1 \to C_x$ a system of special points
$(\gamma(x_0),\dots,\gamma(x_k))$ on the cactus $C_x$.
We build a rooted labelled tree $\Gamma$ almost exactly as in the proof 
of theorem \ref{noncyc}.
There is only a slight difference in the definition of the label of a vertex  associated to a lobe. Let $S_j$ be a lobe of $C_x$ and let 
$H_j \subset S_j$ be the finite subset 
of {\em all} intersection points of $S_j$ with other lobes and of all special points of the form $\gamma(x_i)$ in the lobe. 
We have
a preferred based identification of $S_j$ with $S^1$, but in general 
the base point is not an intersection point or $\gamma(x_0)$ as in the non-cyclic case. 
In the next figure we give an example for $m=3$ and $k=3$. 
The base point of the cactus is $P$ and the base points of the lobes are
$P_1,P_2,P_3$. Only the vertices of the lobes are indicated.
\epsfxsize=250pt
$$\epsfbox{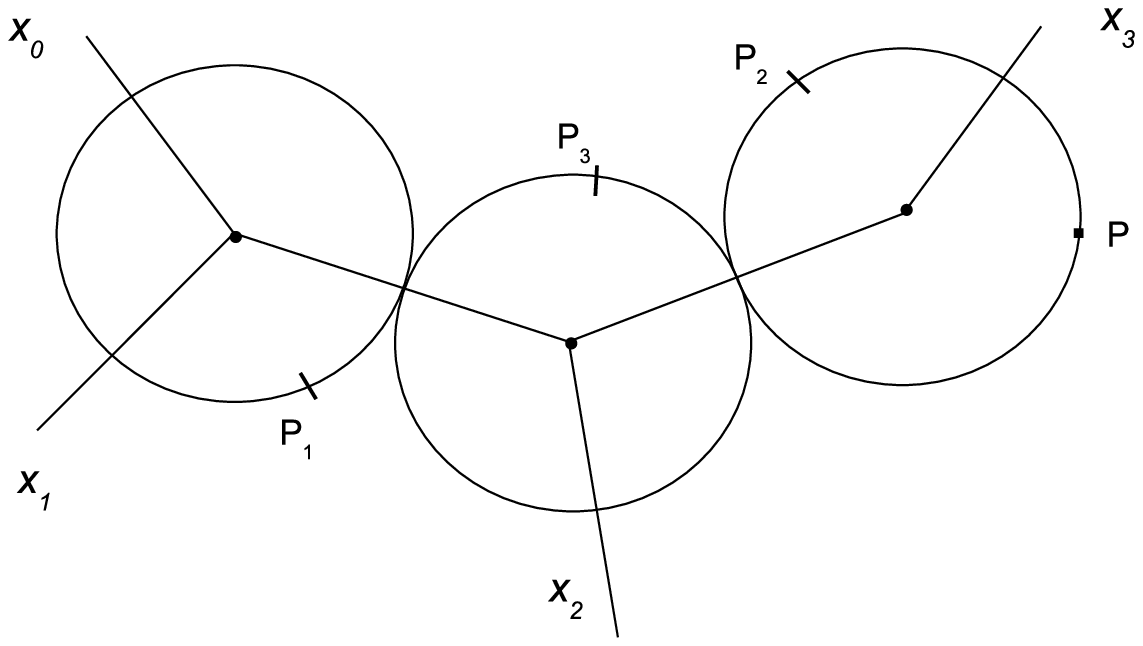}$$ 
Under the identification $S^1 \cong S_j$,  the subset $H_j \subset S_j$ corresponds to a cyclically ordered sequence of points on $S^1$, that is an element $y_j \in \Lambda^{l_j}=\Delta^{l_j} \times S^1$ (actually $y_j$ is in the interior), where
$l_j$ is one less than the cardinality of $H_j$.
The map $a_j^{l_j}:\Lambda^{l_j} \to O(l_j)$ defines the label 
$a_j^{l_j}(y_j) \in O(l_j)$ of the vertex of $S_j$.
We define $a^k(x_0,\dots,x_k) \in O(k)$
as the composition of all labels of the rooted labelled tree $\Gamma$.
The crucial point is that the map $a^k:\Lambda^k \to O(k)$ is equivariant with respect to the action of $\Z_{k+1}$, because $O$ is a cyclic operad.
Moreover
the map $a^*:\Lambda^* \to O^*$ respects cofaces and codegeneracies
because each $a_i^*$ does so.
Finally the structure maps $\theta_*$ make $Tot(O^*)$ into a $fMS$-algebra,
since $O$ is an algebra over the triple $F^+_*$.
\end{proof}

\bex
For a group $G$ the operad $\ud{G}$ is a cyclic operad with
multiplication. The operad with multiplication is in example \ref{topmon}
with $M=G$, and the cyclic group action in example \ref{freeloop}.
These are compatible because $$t_2(m_2)=t_2(1,1)=(1,1)=m_2\, .$$
By theorem \ref{main} \[Tot(\ud{G}^*) \cong \Omega G \simeq \Omega^2 BG\]
has the action of an operad equivalent to the little framed 2-discs.
By proposition \ref{equireco} 
$BG$ is equivalent to a based space with a based $S^1$-action.
We expect the action to be trivial, up to a zig-zag of based equivariant maps that are non-equivariant weak equivalences.

Actually there is a specific based circle action on $BG$, that we
also expect to be equivalent.
A {\em cyclic space} is a {\em contravariant} functor from $\Lambda$ to the category of topological spaces, and in particular it has the underlying 
structure of a simplicial space. The realization of a cyclic space is equipped with a circle action \cite{JJ}.
We recall from \cite{JJ} that the space $LBG$ is equivalent 
(by a $S^1$-equivariant map that is a non-equivariant equivalence)
to the realization of a cyclic space $lbG$ with $(lbG)_n=G^{n+1}$, and \begin{align*}
d_i(x_0,\dots,x_n)&=(x_0,\dots,x_{i-1},x_i x_{i+1},x_{i+2},\dots,x_n)\\
s_i(x_0,\dots,x_n)&=(x_0,\dots,x_{i-1},e,x_{i+1},x_{i+2},\dots,x_n)\\
t_n(x_0,\dots,x_n)&=(x_n,x_0,x_1,\dots,x_{n-1})
\end{align*}
Here $e$ is the unit of $G$.
The realization $|lbG|$ of $lbG$ can be thought of as the space of finite configurations of points on the circle with summable labels in $G$ \cite{Sal}.
There is a cyclic subspace $bG \subset lbG$ that at level $n$
is the space $(bG)_n$  of $n$-tuples $(g_0,\dots,g_n)$ such that the
ordered product satisfies $g_0 \cdots  g_n =1$.
The homeomorphism $(bG)_n \cong G^{n}$ induced forgetting $g_0$ gives
an isomorphism between $bG$ and the classical simplicial construction
$B_*G$ with realization $|B_*G|=BG$.
The category $\Lambda$ is well known to be non-canonically isomorphic to its opposite $\Lambda^{op}$. The choice of an isomorphism gives a bijective correspondence between cyclic spaces and cocyclic spaces.
In particular under the isomorphism described in section 1 of 
\cite{JJ} the cyclic space $bG$ corresponds to the cocyclic space $\omega G$ of example \ref{basedloop}.
\eex

\bex
{\em String topology in the category of correspondences}

The cocyclic space $lX$ does not come from a cyclic operad because 
of course the totalization $LX$ is not a double loop space in general.
However if we work in the category $Cor$ of correspondences this works.
Objects in $Cor$ are topological spaces and a morphisms $X \to Y$ is a pair of maps $X \stackrel{f}{\longleftarrow} Z \stackrel{g}{\longrightarrow} Y$.
The composition is defined via pullbacks. 
When the backward map $f$ is an inclusion we can view a morphism as a 
map defined only partially.
If $f$ is a homeomorphism then we identify the morphism
to the honest map $g \circ f^{-1}:X \to Y$.
The cartesian product induces a monoidal structure on $Cor$ and it makes
sense to consider operads in $Cor$.

It turns out that the collection $lX$ 
forms a cyclic operad with multiplication in the category of correspondences.

The composition map $$\circ_i:(lX)_m \times (lX)_n \to lX_{m+n-1}$$ is defined 
on each pair $((x_0,\dots,x_m),(y_0,\dots,y_n))$ such that
$x_i=y_0$, sending it to 
$$(x_0,\dots,x_{i-1},y_1,\dots,y_n,x_{i+1},\dots,x_m)$$

The operad map from the associative operad $Ass \to lX$ is defined 
at level $n$ by the correspondence 
$* \leftarrow X \stackrel{\Delta}{\rightarrow} X^{n+1}$, where $\Delta$ is the inclusion of the 
thin diagonal. The associated cosimplicial object is the honest
cosimplicial space $lX$.
A small variant of theorem \ref{main} shows that the totalization $LX=Tot(lX)$
admits the action of the operad $fMS$ {\em in the category of		 correspondences}. This is exactly the action constructed by Sasha Voronov \cite{Vor}, at least on the suboperad $Cacti \subset fMS$, and that motivated him to work in the category of correspondences.
\eex

{\em Further directions}

Kaufmann \cite{moduli}
and Kontsevich-Soibelman \cite{KonSoi} proved that
the Hochschild complex of a Frobenius algebra has an action
on the chain level of the operad $\mathcal{M}$ of moduli spaces of Riemann surfaces with boundary. The restriction of this action to the operad of genus 0 surfaces is the cyclic Deligne conjecture. The generalization works because the endomorphism
operad of a Frobenius algebra is a {\em modular operad} \cite{book} with multiplication. Let $O$ be a modular operad in based spaces.
It is natural to ask whether the totalization $Tot(O^*)$ admits
an action of the operad $\mathcal{M}$.

\section{Appendix: comparison with the McClure-Smith construction}

We give the proof of proposition \ref{msd}. More generally we produce for any $m$ a geometric construction of an operad $MS_m$ that is isomorphic to  
the $E_m$ operad  called $\D_m$ by 
McClure-Smith in \cite{MS02}. 
In this section $\D_m$ will denote such operad and not the little $m$-discs operad.
\bde 
Let $MS_m(n)=\F_m(n) \times Mon(I,\partial I)$, where 
$\F_m(n)$ is defined similarly as $\F(n)$ in definition \ref{mani} except that the 4-tuple $(z_1,z_2,z_3,z_4)$ is replaced by a $(m+2)$-tuple $(z_1,\dots,z_{m+2})$ with $z_{2i-1} \in \stackrel{\circ}{I_j}(x)$
and $z_{2i} \in \stackrel{\circ}{I_k}(x)$ for all $i$.
The space $\F_m(n)$ is a regular CW-complex with a cell for each sequence assuming all values between 1 and $n$, with no equal adjacent values, and with no subsequence 
of length $m+1$ of the form $abab\dots$, for $a \neq b$.
The cell is identified to the product
$\prod_{j=1}^n \Delta^{d(j)}$,
where $d(j)$ is one less than the number of occurrences of $j$ 
in the sequence. 
Clearly $\F(n)=\F_2(n)$ and $MS=MS_2$.
 
The direct limits are denoted $MS_\infty=lim_m MS_m$ and
$\F_\infty(n)=lim_m \F_m(n)$. 
There is a symmetric operad structure on $MS_m$ that is
defined as for $MS$ in section \ref{cac} by an embedding onto a suboperad
of $Coend(S^1)$.
\ede
\begin{prop}
There is an isomorphism between the operad $MS_\infty$ and the $E_\infty$ operad
$\D$ in \cite{MS02}. This restricts for each $m$ to an isomorphism between $MS_m$ and the $E_m$ operad denoted $\D_m$ 
in \cite{MS02}. 
\end{prop}
\begin{proof}
By definition $\D(k)=Tot(Y_k^*)$, where $Y_k^*=
\Xi_k(\Delta^*,\dots, \Delta^*)$
 is the cosimplicial space described in section 6 of \cite{MS02}.
Given 
cosimplicial spaces $X_1^*,\dots,X_k^*$, the degree $s$
space $\Xi_k(X_1^*,\dots,X_k^*)^s$ is the colimit 
of a functor $G_1:A_1 \to Top$ to the category of topological spaces. Here $A_1$ is the category whose objects are pairs $(f,h)$ consisting
of a monotone map $f:[t] \to [s]$ and a map $h:[t] \to   \{1,\dots,k\}$,
for some $t$. A morphism from $(f,h)$ to $(f',h')$ is an ordered map
$g:[t] \to [t']$ between the respective domains such that $f' \circ g=f$ and $h' \circ g = h$.   
The functor $G_1$ sends a pair $(f,h)$ to 
$\prod_{1 \leq i \leq k} X^{|f^{-1}(i)|-1}$.

In the special case $X_1^*=\dots=X_k^*=\Delta^*$  corollary 12.3 in \cite{MS02} states that $Y_k^s$  is a regular CW complex with a cell 
for each pair $(f,h)$ as above such that $f$ 
is surjective and for $0 \leq j <t$ either $f(j) \neq f(j+1)$ or $h(j) \neq h(j+1)$.
Such cell has dimension  $t+1-k$ and is  $G_1(f,h)=\prod_{1 \leq i \leq k}  \Delta^{|f^{-1}(i)|-1}$. 
Thus there is a natural identification between the 
cell of $Y_k^0$ indexed by $(f,*)$  and the cell of $\F_\infty(k)$ indexed by the sequence $f(0),\dots,f(t)$.  Here $*:[t] \to [0]$ is the unique map. These identifications are compatible with the inclusion of boundary cells
and give a cellular $\Sigma_k$-equivariant homeomorphism 
$\alpha: \F_\infty(k) \cong Y_k^0$.
In Proposition 12.7 of \cite{MS02} McClure-Smith construct 
an isomorphism of cosimplicial spaces $\omega: Y_k^* \cong \Delta^* \times Y_k^0$.
We describe geometrically the homeomorphism
$$\Phi= \omega^{-1} \circ( \Delta^s \times \alpha): \Delta^s \times \F_\infty(k)  \cong Y_k^s.$$  
A partition $x \in \F_\infty(k)$ is a subdivision of $S^1$ into $k$ $1$-manifolds.  
An element of the simplex 
$z \in \Delta^s$ determines $s$ points on the circle and generically $s$ intervals that we call {\em coarse} intervals, inducing a subdivision of
the 1-manifolds into intervals that we call {\em fine} intervals.
We number by $0,\dots,t$  in the
counterclockwise ordering the fine intervals. 
Let 
$f:[t] \to \{1,\dots,k\}$ be the map sending the index of a fine interval to the index of the 1-manifold containing it. Similarly 
let $h:[t] \to [s]$ be the map sending the index of a fine interval
to the index of the coarse interval containing it. The cell of $\Phi(z,x) \in Y_k^s$
 is indexed by the pair $(f,h)$.
The point $\Phi(z,x)$ of the cell 
$\prod_{1 \leq i \leq k}  \Delta^{|f^{-1}(i)|-1}$
is determined by the 
lengths of the fine intervals: the
$i$-th component corresponds to the lengths of the fine intervals contained in the $i$-th 1-manifold, that sum up to 1 and thus
determine an element of $\Delta^{|f^{-1}(i)|-1}$.
The homeomorphisms $\Phi$ induce a $\Sigma_k$-equivariant homeomorphism 
$$\Psi: \D(k) = Tot(Y_k^*) 
\cong Tot(\Delta^*) \times \F_\infty(k) = MS_\infty(k).$$
Given $f \in \D(a)$ and $g \in \D(b)$, i.e.
cosimplicial maps 

$f^*:\Delta^* \to \Xi_a(\Delta^*,\dots,\Delta^*)=Y_a^*$
and $g^*:\Delta^* \to \Xi_b(\Delta^*,\dots,\Delta^*)=Y^*_b,$
we recall that
the operadic composition $f \circ_i g \in \D(a+b-1)$ is the cosimplicial map
$$\Delta^* \stackrel{f^*}{\rightarrow} Y_a^*  \stackrel{\Xi_a(id,\dots,g^*,\dots,id)}{\longrightarrow}
 \Xi_a(\Delta^*,\dots,\Xi_b(\Delta^*,\dots,\Delta^*),\dots,\Delta^*)
 \stackrel{\Gamma}{\cong} \Xi_{a+b-1}(\Delta^*,\dots,\Delta^*)$$
where $\Gamma$ is the cosimplicial isomorphism described in section
6 of \cite{MS02}. 
Via the geometric correspondence $\Phi$ it is straightforward to check 
that for any $x \in \Delta^s$ 
$$\Psi(f \circ_i g)(x)=[\Psi(f) \circ_i \Psi(g)](x) \in 
\Delta^s \times \F_\infty(a+b-1),$$ and so 
$\Psi:\D \cong MS_\infty$ is an operad isomorphism.

For a fixed $m$ let $Z^k_* \subset Y^k_*$ be the cosimplicial 
subspace in section 13 of \cite{MS02}. 
The $\Sigma_k$-equivariant homeomorphism $\alpha: \F_\infty(k) \cong Y^k_0  $
restricts to a $\Sigma_k$-equivariant homeomorphism $\F_m(k) \cong Z^k_0$, because the cells of $\F_m(k)$ correspond to the cells of $Y^k_0$ whose indexing sequence has {\em complexity} not exceeding $m$,
i.e. the cells of $Z^k_0$.
The splitting $\omega: Y^k_* \cong \Delta^* \times Y_k^0 $ restricts to 
a splitting  $Z^k_* \cong \Delta^* \times Z_k^0$.
The $E_m$-operad $\D_m$ in \cite{MS02} has as
$k$-th space $\D_m(k)=Tot(Z^k_*)$.   
Thus the operad isomorphism $\Psi:\D \cong MS_\infty$ restricts to 
an isomorphism between the suboperad $\D_m \subset \D$  in \cite{MS02} and the suboperad $MS_m \subset MS$. 
\end{proof}

Proposition \ref{msd} follows as a special case for $m=2$.


\begin{thebibliography}{99}

\bibitem{BV} M. Boardman and R. Vogt, Homotopy invariant algebraic structures on topological spaces. Lecture Notes in Mathematics, Vol. 347. Springer-Verlag

\bibitem{Bu} R. Budney, The operad of framed discs is cyclic. 
J. Pure Appl. Algebra 212 (2008), no. 1, 193--196. 

\bibitem{CohVor} R. Cohen and A. Voronov, Notes on string topology,
arXiv:math.AT/0503625

\bibitem{Drin} V. Drinfeld, On the notion of geometric realization, 
Moscow Math. J. 4(2004), no. 3, 619--626.

\bibitem{Ger} M. Gerstenhaber, The cohomology structure of an associative ring.  Ann. of Math. (2)  78  1963 267--288.

\bibitem{getzler} E. Getzler, Batalin-Vilkovisky algebras and two-dimensional topological field theories.  Comm. Math. Phys. 159 (1994), no. 2, 265--285.

\bibitem{gk} E. Getzler and M. Kapranov, Cyclic operads and cyclic homology,  Geometry, topology, \& physics,  167--201, Conf. Proc. Lecture Notes Geom. Topology, IV, Int. Press, Cambridge, MA, 1995.

\bibitem{Hu} P. Hu, The Hochschild cohomology of a Poincar\'e algebra,
arXiv:0707.4118

\bibitem{JJ} J. Jones, Cyclic homology and equivariant homology.  Invent. Math. 87 (1987), no. 2, 403-423. 

\bibitem{Ka} R. Kaufmann, On several varieties of cacti and their relations. Algebraic and Geometric Topology 5 (2005), 237-300.

\bibitem{KaDel} R. Kaufmann, A proof of a cyclic version of Deligne's conjecture via cacti, arXiv:math.QA/0403340

\bibitem{moduli} R. Kaufmann, Moduli space actions on the Hochschild co-chains of a Frobenius algebra I-II, arXiv:math.AT/0606064-0606065

\bibitem{KonSoi} M. Kontsevich and Y. Soibelman, 
Notes on A-infinity algebras, A-infinity categories and non-commutative geometry. I, arXiv:math.RA/0606241

\bibitem{book} M. Martin, S. Shnider and J. Stasheff,
Operads in algebra, topology and physics. 
Mathematical Surveys and Monographs, 96. American Mathematical Society.

\bibitem{May} P. May, 	The geometry of iterated loop spaces. Lectures Notes in Mathematics, Vol. 271.

\bibitem{MS99} J. McClure and J. Smith, A solution of Deligne's Hochschild cohomology conjecture. Recent progress in homotopy theory (Baltimore, MD, 2000),  153--193, Contemp. Math., 293, Amer. Math. Soc.

\bibitem{MS02} J. McClure and J. Smith, Cosimplicial objects and little $n$-cubes. I.  Amer. J. Math.  126  (2004),  no. 5, 1109--1153.

\bibitem{MS04} J. McClure and J. Smith, Operads and cosimplicial objects: an introduction.  Axiomatic, enriched and motivic homotopy theory,  133--171, NATO Sci. Ser. II Math. Phys. Chem., 131, Kluwer Acad. Publ., Dordrecht, 2004. 

\bibitem{Men} L. Menichi, Batalin-Vilkovisky algebras and cyclic cohomology of Hopf algebras.  $K$-Theory 32 (2004), no. 3, 231--251.

\bibitem{Sal} P. Salvatore, Configuration spaces with summable labels.  Cohomological methods in homotopy theory 375--395, Progr. Math., 196, Birkh\"auser, Basel, 2001. 

\bibitem{imrn} P. Salvatore, Knots, operads, and double loop spaces.
Int. Math. Res. Not. 2006, Art. ID 13628. 

\bibitem{SW} P. Salvatore and N. Wahl, Framed discs operads and
Batalin-Vilkovisky algebras, Q.J.Math. 54 (2003), 213-231.

\bibitem{Sin} D. Sinha, Operads and knot spaces, J. Amer. Math. Soc. 19 (2006) 461-486. 

\bibitem{Vor} A. Voronov, Notes on universal algebra. Graphs and patterns in mathematics and theoretical physics,  81--103,
Proc. Sympos. Pure Math., 73, Amer. Math. Soc., Providence, RI, 2005. 


\end{thebibliography}
\end{document}